\documentstyle[12pt]{article}

\begin{document}
\newtheorem{guess}{Proposition }[section]
\newtheorem{theorem}[guess]{Theorem}
\newtheorem{lemma}[guess]{Lemma}
\newtheorem{corollary}[guess]{Corollary}
\def \ja {\vrule height 3mm width 3mm}
\def\inftwo#1#2#3{\smash{\mathop{#1}\limits_{#2\atop{#3}}}}

\centerline{\Large \bf Chern Classes  of Bundles on Rational Surfaces}
\vspace{.2 in}
\centerline{\Large Elizabeth  Gasparim }
\centerline{\small Departamento de Matem\'atica, Universidade Federal de Pernambuco} 
\centerline {\small Cidade Universit\'aria, Recife, PE, BRASIL, 50670-901}
\centerline{\small gasparim@dmat.ufpe.br}

\begin{abstract}
Consider the blow up 
 $\pi: \widetilde{X} \rightarrow X$  of a rational
surface
$X$ at a point.
Let  $\widetilde{V}$ be a holomorphic
 bundle over $\widetilde{X}$
whose restriction to the exceptional divisor 
is ${\cal O}(j) \oplus {\cal O}(-j)$ 
 and define $V =(\pi_*\widetilde{V})^{\vee \vee}.$
Friedman and Morgan gave the following 
bounds for the second Chern classes
$j \leq c_2(\widetilde{V}) - c_2(V) \leq j^2.$
We show that these bounds are sharp.
\end{abstract}

\section{Introduction.}

In this paper our basic setting will be the following.
  $X$ will be a rational surface, $\pi : \widetilde {X} \rightarrow X$
 the blow-up of $X$  at point $x  \in X$ and  
$\ell$ the 
exceptional divisor.
 $\widetilde {V}$ will be a rank two 
bundle  over the surface
$\widetilde X$ satisfying
$det\, \widetilde {V} \simeq {\cal O}_ {\widetilde X}$ 
and $\widetilde{V}|_{\ell}  \simeq {\cal O}(j) \oplus
 {\cal O}(-j), \, j \geq 0.$
Let $V = \pi_* {\widetilde{V}}^{\vee \vee}.$
Friedman and Morgan [2] gave the following estimate 
for the second Chern classes
$j \leq  c_2(\widetilde{V}) - c_2(V) \leq j^2.$
We show that
these
bounds are sharp by
 giving an algebraic procedure to 
calculate 
$c_2(\widetilde{V}) - c_2(V)$ directly from the transition 
matrix defining $\widetilde{V}$ 
in a neighborhood of $\ell.$

Since  $\widetilde{V}|_{\widetilde{X} - {\ell}} =
\pi^* V|_{X - p} $            
it is natural to localize our  
study of $\widetilde{V}$
 to a neighborhood of the exceptional divisor.
If $x$ is the blown-up point, we
choose an open set $U \ni x$
that is biholomorphic to ${\bf C}^2.$ 
Then $\pi^* (U)$  gives
 a neighborhood $N_{\ell} \sim \widetilde{{\bf C}^2}.$
We remark that since $\widetilde{X}$ is rational 
we can only guaranty the existence of a 
neighborhood of $x$  that is biholomorphic to ${\bf C}^2$
minus a finite number of points. However, all 
our calculations use holomorphic functions, which always extend over    
these points since we are in dimension 2.
Therefore me may from the start assume 
without loss of generality that $N(\ell) = 
\widetilde{{\bf C}^2}.$
We  use the following results:

 \begin{theorem}:([5], Thm. 2.1)   Let $E$ be a
 holomorphic bundle on
 $ \widetilde{\bf C}^2  $ 
with
$E_{\ell} \simeq {\cal O}(j) \oplus  {\cal O}(-j). $
Then  $E$  has a transition matrix
of the form
$\left(\matrix {z^j & p \cr 0 &  z^{-j} \cr }\right)$
from $U$ to $V,$  where 
$p = \sum_{i = 1}^{2j-2} \sum_{l = i-j+1}^{j-1}p_{il}z^lu^i.$
\end{theorem}

\begin{corollary}([6], Cor. 4.1) \label{topzer}Every
holomorphic rank two vector bundle over $\widetilde{X}$
with vanishing first Chern class
is topologically  determined
 by a triple $(V,j,p)$
where $V$ is a rank two holomorphic  bundle on $X$ with
vanishing first Chern class, $j$ is a non-negative
integer, and  $p$ is a polynomial.
\end{corollary}

It follows that for each polynomial $p$ ( in 3 variables 
$z,$ $z^{-1},$ and $u$ ) there is a canonical construction that 
assings a bundle $\widetilde{V}$ over $\widetilde{X}$ whose topological
type is determined by $p.$
We remark that N. Buchdahl [1] has 
shown that these bounds are sharp in a more 
general setting by entirely different methods.

Following Friedman and Morgan [2, p. 302], we 
 define the sheaf  $Q$ by the
exact sequence
\begin{equation}\label{seq}
0 \rightarrow \pi_* \widetilde {V} \rightarrow
V \rightarrow Q \rightarrow 0. \label{eq- seq} 
\end{equation}
where $Q$ is a sheaf supported only at $x.$
 From the exact sequence (\ref{seq})
it
follows immediately that
 $ c_2\pi_*(\widetilde{V}) - c_2(V) = l(Q),$
where $l$ stands for length.
An application of
Grothendieck-Riemann-Roch (cf [2]) gives
that $c_2( \widetilde{V}) - c_2(V) = l(Q) +
l(R^1 \pi_* \widetilde{V}).$
 Since $Q$ and $R^1\pi_* \widetilde {V}$ are supported at $x$
we can compute their lengths by looking  at $\widetilde {V}|_{N_{\ell}}$
and at $\pi$ as the blow-up map
$\pi: \widetilde{\bf C}^2 \rightarrow {\bf C}^2.$
The lengths of
$Q$ and $R^1\pi_* \widetilde{V}$
can be explicitly computed using the 
Theorem on Formal Functions of Grothendieck (see[3]).

\section{ Calculation of Chern classes.}

\subsection{The upper bound occurs for $p =0.$}

\begin{guess}\label{j^2}: If the bundle $\widetilde{V}$
splits on a neighborhood of the exceptional divisor,
then  $c_2(\widetilde{V}) - c_2(V) = j^2.$
\end{guess}

\noindent {\bf Proof}: 
We  show that  for
 $\widetilde{V}|_{N(\ell)} = {\cal O}(j) \oplus {\cal O}(-j),$
we have
$l(R^1\pi_* \widetilde{V}) = j(j+1)/2$ 
and $l(Q) = j(j-1)/2,$ 
which we state as  lemmas.
It follows that 
$c_2(\widetilde{V}) - c_2(V)  = 
l(Q) +  l(R^1\pi_* \widetilde{V}) = j(j+1)/2 + j(j-1)/2 = j^2.$\hfill\ja

We now prove the lemmas we just used.

\begin{lemma}\label{l(Q)}
 If the bundle $\widetilde{V}$ splits on a neighborhood 
of the exceptional divisor,  then 
$l(Q) = j(j+1)/2.$
\end{lemma}

\noindent{\bf Proof}:
 Since  the length of
$Q$ equals the dimension
of $Q_x^{\wedge}$ as a $k(x)$-vector space, we
need to study the map $ (\pi_*  \widetilde {V}^{\wedge}_x)
 \rightarrow V^{\wedge}_x$
and compute the dimension of the  cokernel as  a $k(x)$-vector space.
But as $V^{\wedge}_x = (\pi_* \widetilde{V}^{\wedge}_x)^{\vee \vee},$
we need to compute the ${\cal O}^{\wedge}_x$-module structure
on $M = (\pi_* \widetilde{V_x})^{\wedge}$ and study the natural map
$M \hookrightarrow M^{\vee \vee}$ of ${\cal O}^{\wedge}_x$-modules.
By the Formal Functions Theorem
$$M \simeq   \lim_{\longleftarrow}H^0(\ell_n, \widetilde{V}|{\ell_n})$$
as   ${\cal O}^{\wedge}_x ( \simeq
{\bf C}[[x,y]]$)- modules, where
 $\ell_n \simeq N_{\ell} \times {\cal O}_x/m^{n+1}_x$
is the n-th infinitesimal neighborhood of $\ell.$
To do this we will have to isolate the action of
$\displaystyle{{\bf C}[[x,y]]} \over{ (x,y) ^{n+1}}$  on
$H^0(\ell_n, \widetilde{V}|{\ell_n}).$

We first write the blow-up of ${\bf C}^2$ with two charts 
$U \sim  V \sim {\bf C}^2 $ with $(z,u) \mapsto (z^{-1}, zu)$
in $U \cap V.$ Then the blow-up map
$\pi : \widetilde{{\bf C}^2} \rightarrow {\bf C}^2$
is given on the $U$ chart by 
$ (x,y) = \pi(z,u) = (u,zu).$
We give the natural action of $x$  and $y$ on this space; that is, $x$
acts by multiplication by $u$ and 
$y$ acts by multiplication by $zu.$ This yields 
$M = {\bf C}[[x,y]]<\alpha,\beta_0,\beta_1,...,\beta_j>,$
where 

$$\alpha = \left(\matrix{ 1 \cr 0}\right),\,\, 
\beta_0 = \left(\matrix{ 0 \cr 1}\right),\,\,
\beta_1 = \left(\matrix{ 0 \cr z}\right),\,\, \cdots,\,\,
\beta_j = \left(\matrix{ 0 \cr z^j}\right).$$
With relations:

$$\left\{\matrix {x\beta_1 - y \beta_0 \cr
x \beta_2 - y \beta_1 \cr
\vdots \cr
x\beta_j - y \beta_{j-1}}\right.,
\,\, \,
\left\{\matrix{x^2\beta_2 - y^2 \beta_0 \cr
\vdots \cr
x^2\beta_j - y^2 \beta_{j-2}}\right.,
\,\, \,\cdots \,\,\,,
\,\,\, \left\{\matrix{ x^j \beta_j - y^j \beta_0.}\right.
$$
All together there are $j(j+1)/2$ relations.
Now writing the generators of $M^{\vee}$ and
$M^{\vee \vee}$ we see that $coker(M \hookrightarrow M^{\vee \vee})$
is a $j(j+1)/2$ dimensional vector space over ${\bf C} \simeq k(x),$ 
hence $l(Q) = j(j+1)/2.$\hfill\ja
	
\begin{lemma}\label{oplus}
 If we have a split extension on a neighborhood 
of the exceptional divisor, that is, when $\widetilde{V}|_{N(\ell)} \simeq
{\cal O}(j) \oplus {\cal O}(-j),$ then 
$l(R^1\pi_* \widetilde{V}) = j(j-1)/2.$
\end{lemma}

\noindent {\bf Proof}:
Here again we follow the same method 
using the Theorem on Formal Functions and
compute that 
$$h^1(\ell_n, \widetilde{V}|_{\ell_n}) = 
\left\{\begin{array}{ll} 
(j-1)\,\, &  for \,\, n = 0 \\
(j-1)+(j-2)+ \cdots + (j-n)\,\, &  for \,\, 1\leq n \leq j-1 \\
j(j-1)/2 \,\, & for \,\, n \geq j-1. \end{array} \right.$$

Also $ Ker \left(H^1(\ell_n, \widetilde{V}|_{\ell_n}) \rightarrow
H^1(\ell_{n-1},\widetilde{V}|_{\ell_{n-1}})\right)$ has dimension
$(j-n)$ for $2 \leq n \leq j-1.$ 
Computing the inverse limit (cf. Lang [4]) this implies that 
$$ l(R^1\pi_* \widetilde{V}) =  dim\,\, \inftwo{\lim}{\longleftarrow}{n}\,\,
 H^1(\ell_n, \widetilde{V}|_{\ell_n}) = 
1+2+ \cdots + (j-1) = j(j-1)/2.$$
This can also be verified by  invoking 
the  short exact sequence on p.388 of Hartshorne [3].\hfill\ja

\vspace{5 mm}

\noindent {\bf Remark}: We can also proof Proposition \ref{j^2}
in a simpler way, by 
explicitly constructing a generic section  of 
$\widetilde{V}$ and
counting its zeros.

\subsection{The lower bound occurs for $p =u.$}
\begin{theorem} If $\widetilde{V}$ is a bundle corresponding to the 
triple $(V,j,u)$ (according to Cor. 1.2), then 
$c_2(\widetilde{V})- c_2(V) = j.$
\end{theorem}

\noindent {\bf Proof}: 
We  show that  
$l(R^1\pi_* \widetilde{V}) = j-1$ 
and $l(Q) = 1,$
which we state as  lemmas.
It follows that 
$c_2(\widetilde{V}) - c_2(V)  = (j-1) +1 = j.$\hfill\ja

\begin{lemma} If $\widetilde{V} $ is given by $(V,j,u)$ 
then $l(Q) = 1.$
\end{lemma}

\noindent {\bf Proof}:
The bundle  $\widetilde{V}$ is given over $N(\ell)$
(according to Theorem 2.1) by the transition matrix
$\left(\matrix{ z^j & u \cr 0 & z^{-j}}\right).$
Here  calculations are similar to  the ones in the
proof of Lemma 2.2.
We set $M = (\pi_* \widetilde{V_x})^{\wedge}$ and study the natural map 
$\rho: M \hookrightarrow M^{\vee \vee}$ of ${\cal O}^{\wedge}_x$-modules.
The value of $l(Q)$ is the dimension of cokernel of $\rho.$

The ${\cal O}_x$ - module structure of $M$ is completely determined 
by the structure of the sections of the bundle $V'.$
Writing  sections of $V'$ in the form 
$\left( \matrix{ a \cr b}\right) = 
\left( \matrix{ \sum a_{ik}z^ku^i \cr \sum b_{ik}z^k u^i }\right)$
it is simple to see which restrictions are imposed in the coefficients
$a_{ik}$ and $b_{ik}.$ Some terms of the  general form for the
sections are 
$$\left( \matrix{ a \cr b}\right) = 
b_{00}\left( \matrix{ 0 \cr 1}\right)
 + b_{01}\left( \matrix{ 0 \cr z}\right)
+b_{0j}\left( \matrix{ -u \cr z^j }\right) +
\cdots + 
a_{j0}\left( \matrix{ u^j \cr 0 }\right)+ \cdots .$$
In fact, one verifies that these terms are enough to generate all
 the sections. It then follows that
$M = <\beta_{00}, \beta_{01},\beta_{0j},\alpha{jo}>/R$
where 
$\beta_{00} = \left( \matrix{ 0 \cr 1}\right),
 \beta_{01} = \left( \matrix{ 0 \cr z}\right),
 \beta_{0j} = \left( \matrix{ -u \cr z^j }\right),
 \alpha_{j0} = \left( \matrix{ u^j \cr 0 }\right)$ 
and $R$ is the set of relations
$$\left\{ \begin{array} {l} 
x \, \beta_{01} - y \,  \beta_{00} = 0 \cr
\alpha_{j0} + x^{j-1} \beta_{0j} - y^{j-1}\beta_{01} = 0\cr
\end{array}\right..$$
Using the second relation, one eliminates
$\alpha_{j0}$ from the set of generators and 
gets a simpler presentation 
 $M \simeq  <\beta_{00}, \beta_{01},\beta_{0j}>/R'$
where $R'$ now has the single relation 
$x \, \beta_{01} - y \,  \beta_{00} = 0.$
It is now a matter of simple algebra to find that 
$M^{\vee} = < a, b > $ is free on two generators, where
$a = \left\{\begin{array} {l} \beta_{00} \mapsto x \cr
                              \beta_{01} \mapsto y \cr 
                              \beta_{0j} \mapsto 0 \cr
     \end{array}\right.$ and 
$b = \left\{\begin{array} {l} \beta_{00} \mapsto 0 \cr
                              \beta_{01} \mapsto 0 \cr 
                              \beta_{0j} \mapsto 1 \cr
     \end{array}\right..$
Then naturally $M^{\vee\vee} = <  a^* , b^*>$ 
is generated by the dual basis, namely 
$a^* = \left\{\begin{array} {l} a \mapsto 1 \cr
                                b  \mapsto 0 \cr 
     \end{array}\right.$ and 
$b^* = \left\{\begin{array} {l} a \mapsto 0 \cr
                                b  \mapsto 1 \cr
     \end{array}\right..$
The map $\rho$ is given by evaluation and we have 
$im \,\rho = < x\, a^*, y \, a^*, b^*>$ and therefore the 
cokernel is $coker \, \rho = < \overline{a^*}>$
and $l(Q) = dim \, coker \, \rho = 1.$ \hfill\ja

\begin{lemma} If $\widetilde{V} $ is given by $(V,j,u)$ 
then $l(R^1\pi_*\widetilde{V}) = j-1.$
\end{lemma}

\noindent {\bf Proof}: 
We claim that $H^1(\ell_n, \widetilde{V}|_{\ell_n}) $
is generated by the  1-cocycles  
$\left(\matrix{z^k \cr 0}\right)$ for $ -j \leq k \leq -1$ 
and that the maps 
$H^1(\ell_n, \widetilde{V}|_{\ell_n}) \rightarrow 
H^1(\ell_{n-1}, \widetilde{V}|_{\ell_{n-1}}) $ are the identity.
Hence $l(R^1\pi_*\widetilde{V}) = dim\,\, \inftwo{\lim}{\longleftarrow}{n}\,\,
 H^1(\ell_n, \widetilde{V}|_{\ell_n}) =   j-1.$
In fact, if $T$ is the transition matrix 
for $\widetilde{V}|_{N(\ell)}$ then the equality 
$$B = \sum_{i = 0}^\infty\sum_{k = -\infty}^\infty\left(\matrix{ 0 
 \cr b_{ik} z^ku^i}\right) = 
\sum_{i = 0}^\infty\sum_{k =0}^\infty \left(\matrix{ 0 
\cr b_{ik}z^ku^i}\right) + 
T^{-1} \sum_{i = 0}^\infty\sum_{k = -\infty}^{-1} \left(\matrix{ 
b_{ik}z^ku^{i+1}  \cr 
b_{ik}z^{k-j}u^i } \right)$$
 shows that $B$ is a coboundary, since the first term on the r.h.s. is 
holomorphic in $U$ and the last term of the r.h.s. is holomorphic in $V.$
As a consequence every 1-cocycle has a representative of the form
$\alpha = \sum_{i = 0}^\infty\sum_{k = -\infty}^\infty \left(\matrix{ a_{ik}z^ku^i \cr 0  }\right).$
Analogously 
$ A = 
\sum_{i = 0}^\infty\sum_{k =0}^\infty \left(\matrix{ a_{ik} z^ku^i
 \cr  0 }\right) + 
T^{-1} \sum_{i = 0}^\infty\sum_{k = -\infty}^{-1} \left(\matrix{ 
a_{ik}z^ku^i \cr 
 0 } \right)$ is a coboundary.  Therefore the only terms that give
nonzero cohomology classes  in $\alpha$ are the ones with indices 
$k$ for $-j \leq k \leq -1$ and  the claim follows. \hfill\ja

Acknowledgements: I would like to thank N. Buchdahl and 
E. Ballico for pointing out a mistake in the original. 
Special thanks go to R. Hartshorne for useful comments.

\end{document}